\documentclass[11pt,a4paper]{article}

\usepackage{epsf,epsfig,amsfonts,amsgen,amsmath,amstext,amsbsy,amsopn,amsthm
}
\usepackage{ebezier,eepic}
\usepackage{color}

\newtheorem{theorem}{Theorem}

\newtheorem{proposition}{Proposition}
\newtheorem{lemma}{Lemma}
\newtheorem{false statement}{False statement}

\theoremstyle{definition}

\newtheorem{corollary}[theorem]{Corollary}

\renewcommand{\mod}{{\rm mod}}

\newcounter{mathitem}
  {\begin{list}{{$(\roman{mathitem})$}}{
   \setcounter{mathitem}{0}
   \usecounter{mathitem}
   \setlength{\topsep}{0pt plus 2pt minus 0pt}
   \setlength{\parskip}{0pt plus 2pt minus 0pt}
   \setlength{\partopsep}{0pt plus 2pt minus 0pt}
   \setlength{\parsep}{0pt plus 2pt minus 0pt}
   \setlength{\leftmargin}{35pt}
   \setlength{\itemsep}{0pt plus 2pt minus 0pt}}}
  {\end{list}}

\begin{document}

\title{\bf\Large Minimal sufficient sets of colors and minimum number of colors \footnote {~Mathematics Subject Classification 2010: 57M25 57M27}}

\date{}

\author{Jun Ge, Xian'an Jin, Louis H. Kauffman, Pedro Lopes and Lianzhu Zhang \footnote {~Corresponding author.}}
\maketitle
\begin{abstract}

In this paper we first investigate minimal sufficient sets of colors for $p=11$
and $13$. For odd prime $p$ and any $p$-colorable link $L$ with $\det L\neq 0$,
we give alternative proofs of $\text{mincol}_p L\geq 5$ for $p\geq 11$ and
$\text{mincol}_p L\geq 6$ for $p\geq 17$. We elaborate on equivalence classes of sets of distinct colors (on a given modulus) and prove that there are two such classes of five colors modulo $11$, and only one such class of five colors modulo 13.
Finally, we give a positive answer to a question raised by Nakamura, Nakanishi,
and Satoh concerning an inequality involving crossing numbers. We show
it is an equality only for the trefoil and for the figure-eight knots.

\medskip
\noindent {\bf Keywords:} Link colorings; Minimum number of colors; Equivalence classes of colorings; minimal sufficient sets of colors.
\smallskip
\end{abstract}

\section{Introduction}

\baselineskip=18pt

A Fox $m$-coloring \cite{Fox} is an assignment of elements from $\{0, 1, \ldots, m-1\}$ to
the arcs of a link diagram such that at each crossing twice the integer assigned
to the over-arc equals to the sum of the integers assigned to the two under-arcs
mod $m$. For each link diagram and each modulus $m>1$, there are always $m$
trivial colorings, namely by assigning the same integer mod $m$ to every arc of
the diagram. A coloring with at least two distinct colors (i.e., two distinct integers mod $m$
assigned to two arcs) is called a non-trivial coloring. It is easy to check
that if one diagram of a link has a non-trivial $m$-coloring, then each diagram
of that link has a non-trivial $m$-coloring. A link is called $m$-colorable if
it admits a diagram with non-trivial $m$-colorings. The following well-known theorem presents a
criterion for checking if a given link is $m$-colorable.

\begin{theorem}[\cite{LM}]\label{determinant}
A link $L$ is $m$-colorable if and only if the determinant of $L$ ($\det L$)
and $m$ are not relatively prime.
\end{theorem}

The following definition was introduced by Harary and Kauffman in \cite{HK}.

\medskip

DEFINITION 1. ~Given an integer $m$ greater than $1$. Let $L$ be a link admitting non-trivial
$m$-colorings. Let $D$ a diagram of $L$, let $n_{m, D}$ be the minimum
number of colors mod $m$ it takes to construct a non-trivial $m$-coloring on $D$.
Set
$$\text{mincol}_m L \doteq \min \{n_{m, D} \mid D\ \text{is a diagram of}\ L \}.$$
We call $\text{mincol}_m L$ the minimum number of colors of $L$, mod $m$.

\medskip

We call any non-trivial $m$-coloring of $L$ using $\text{mincol}_m L$ colors a
minimal $m$-coloring of $L$.

\bigbreak

In this article we give a shorter proof of the following Theorem.

\begin{theorem}[\cite{LM}]
Let $p$ be a prime greater than $7$, let $L$ be a link with $p\mid \det L \neq 0$. Then $mincol_p\, L \geq 5$.
\end{theorem}

We also prove the following fact.

\begin{theorem}
Let $p$ be a prime greater than $13$, let $L$ be a link with $p\mid \det L \neq 0$. Then $mincol_p\, L \geq 6$.
\end{theorem}

We give these proofs because we think they are both insightful and instructive. In fact some of our methods lead to new research problems which we will detail later in this paper.

Nakamura, Nakanishi, and Satoh proved a more general theorem in \cite{NNS} for knots, which we were unaware of at the time of the writing of our article using a completely different approach.

\begin{theorem}\label{NNS}
Let $p$ be an odd prime. Any $p$-colorable knot $K$ satisfies
\begin{equation}
\text{mincol}_p (K)\geq \lfloor \log_2 p \rfloor+2
\label{inequality}
\end{equation}
where $\lfloor x \rfloor$ is the largest integer less than or equal to $x$.
\end{theorem}

It is worth pointing out that Nakamura et al.'s proof for Theorem \ref{NNS} can not be
naturally extended to $p$-colorable links with non-zero determinant. We would like to
understand whether links with non-zero determinant also have such a good lower bound or not.

Let $c(K)$ denote the crossing number of $K$. Since $c(K)\geq \text{mincol}_p (K)$, any
$p$-colorable knot $K$ satisfies $c(K)\geq \lfloor \log_2 p \rfloor +2$. In \cite{NNS},
Nakamura, Nakanishi, and Satoh ask if the equality only holds for the trefoil knot ($p=3$) and
the figure-eight knot ($p=5$). We will give a positive answer to this question for classical knots
in Section 5, Theorem \ref{equality}.

We also define an equivalence relation among sets of colors over a given modulus, see Definition 5, and count the number of such equivalence classes for the least $mincol_p$ for $p=11$ or $13$. This equivalence relation among sets of colors unveils somewhat more the fascinating topic of Fox colorings.

\section{Definitions and Preliminaries}

REMARK 1. All split links (which have zero determinant) can be non-trivially colored with
two colors in any modulus. That is to say, for any split link and any modulus,
the minimun number of colors is $2$. Moreover, links with zero determinant are colorable over any modulus.
On the other hand, the fact that a link is colorable
only over a finite number of moduli introduces a certain regularity that allows
to prove Theorems like Theorem \ref{mincol} below.
For this reason all links referred to in this
paper are links with non-zero determinant.

\medskip

Clearly, $\text{mincol}_m L$ is an invariant of link $L$. For small primes,
we have certain results for this invariant, listed as Theorem \ref{mincol}
\cite{KL1, Satoh, Oshiro, Saito, LM, CJZ}. However, the minimum number of
colors is very difficult to compute in general, even for torus knots
$T(2, n)$ \cite{KL1, KL2}.

\begin{theorem}\label{mincol}
Let $L$ be a link with non-zero determinant.

(1) If $2\mid \det L$, then $\text{mincol}_2 L=2$.

(2) If $3\mid \det L$, then $\text{mincol}_3 L=3$.

(3) If $5\mid \det L$, then $\text{mincol}_5 L=4$.

(4) If $7\mid \det L$, then $\text{mincol}_7 L=4$.

(5) If $p$ is a prime greater than $7$ and $p\mid \det L$, then $\text{mincol}_p L\geq 5$.

(6) If $11\mid \det L$, then $5\leq \text{mincol}_{11} L\leq 6$.
\end{theorem}

More precisely, Satoh \cite{Satoh} proved that any 5-colorable link
with non-zero determinant can be colored by $\{1, 2, 3, 4\}$. Oshiro
\cite{Oshiro} proved that any 7-colorable link with non-zero determinant
can be colored by \{0, 1, 2, 4\}. And Cheng et al. \cite{CJZ} proved that
any 11-colorable link with non-zero determinant can be colored by 5 or 6
colors in \{0, 1, 4, 6, 7, 8\}. These three papers used similar techniques
developed by Satoh.

For $p=11$, the fourth named author of this article found an interesting behavior recently (\cite{Lopes}). In order
to describe this discovery, we list some terminology introduced in \cite{Lopes} first.

\medskip

DEFINITION 2. ~Let $m$ be a positive integer greater than $1$ and let $L$ be a link
admitting non-trivial $m$-colorings.

$\bullet$ An $m$-sufficient set of colors (for $L$) is a set of integers
mod $m$ such that a non-trivial $m$-coloring can be realized on a diagram
of $L$ with colors from this set.

$\bullet$ An $m$-minimal sufficient set of colors (for $L$) is an $m$-sufficient
set of colors (for $L$) whose cardinality is $\text{mincol}_m L$.

\medskip

\medskip

DEFINITION 3. ~Suppose $m$ is a positive integer greater than $1$ such that for any two links,
$L$ and $L'$, admitting non-trivial $m$-colorings, $\text{mincol}_m L = \text{mincol}_m L'$;
in this case, let $\text{mincol} (m)$ denote this common minimum number of colors,
$\text{mincol} (m) = \text{mincol}_m L = \text{mincol}_m L'$.

If there is a set $S$ of $\text{mincol} (m)$ elements from $\mathbb{Z}_m$ such
that for any $m$-colorable link $L''$, there exists a non-trivial $m$-coloring
with colors from this set, then we call $S$ a common $m$-minimal sufficient set
of colors.

If there is a set $S'$ of $\text{mincol} (m)$ elements from $\mathbb{Z}_m$ such
that any diagram supporting minimal $m$-colorings can be colored with colors from
this set, then we call $S'$ a universal $m$-minimal sufficient set of colors.

\medskip

Combining results proved in \cite{KL1, Satoh, Oshiro, Saito, LM, Lopes}, it is proved in \cite{Lopes}
that for each prime $p< 11$, there is a universal $p$-minimal sufficient set of colors.
Moreover, for any of these primes, any common $p$-minimal sufficient set of colors is a
universal $p$-minimal sufficient set of colors.

Changes happen at $p=11$. It is also proved in \cite{Lopes} that

\begin{theorem}\label{Lopes}
There is no universal 11-minimal sufficient set of colors.
\end{theorem}

We remark it remains unknown whether there is a common 11-minimal sufficient set of colors
or not.

To make our results easier to understand, we now introduce more terminology and
results already in the literature.

\medskip

DEFINITION 4. ~Given a positive integer $m$, we define an $m$-coloring automorphism to be a
permutation $f$, of the set $\mathbb{Z}_m$, such that
$$f(a\ast b)=f(a)\ast f(b)$$
for all $a, b\in \mathbb{Z}_m$, with $x\ast y=2y-x\ (\text{mod}\ m)$,
for every $x, y\in \mathbb{Z}_m$.

\medskip

\begin{proposition}
$m$-coloring automorphisms preserve $m$-colorings.
\end{proposition}

\begin{proposition} \cite{EMR}
Given a positive integer $m$, each $m$-coloring automorphism is given by:
$$f_{\lambda, \mu}(x)=\lambda x+\mu$$
with $\mu \in \mathbb{Z}_m$ and $\lambda \in \mathbb{Z}_m^\ast$, the set of
units of $\mathbb{Z}_m$. Furthermore, the set of $m$-coloring automorphisms
equipped with composition of functions constitutes a group.
\end{proposition}

For a permutation $f$ and a set $\{a_i\}_{i\in T}$, we denote $f(\{a_i\}_{i\in T})=\{f(a_i)\}_{i\in T}$.
Let $D$ be an $m$-colorable diagram with arcs $x_1, \ldots, x_n$. Suppose
$C$ is an $m$-coloring of $D$ with arc $x_i$ colored by $a_i$. Let $f$ be
an $m$-coloring automorphism. Similarly, we denote $f(C)$ the $m$-coloring
of $D$ such that arc $x_i$ is colored by $f(a_i)$.

Now it is natural for us to define an equivalence relation among color sets.

\medskip

DEFINITION 5. ~Given a positive integer $m$. Two color sets $S_1$ and $S_2$ are called
equivalent mod $m$ if and only if there exists an $m$-coloring automorphism
$f_{\lambda, \mu}(x)=\lambda x+\mu$ (mod $m$), such that $f(S_1)=S_2$.

\medskip

Since coloring automorphisms preserve colorings, if a link $L$ can be colored
with colors in a color set $S$, then it can be colored with colors in any color
set $S'$ equivalent to $S$. Furthermore, if a color set $S$ is an $m$-sufficient
set ($m$-minimal sufficient set) for $L$, then any color set $S'$ equivalent to
$S$ is an $m$-sufficient set ($m$-minimal sufficient set) for $L$.

\begin{lemma}\label{012}
Let $L$ be a non-split link. Let $p\geq 3$ be a prime number. If a diagram
$D$ of $L$ admits a non-trivial $p$-coloring $C$, then there exist a $p$-coloring
automorphism $f$ such that $f(C)$ is a non-trivial $p$-coloring of $D$ containing colors
0, 1 and 2.
\end{lemma}

\medskip

Proof. ~Since $L$ is a non-split link and $p$ is an odd prime, $D$ must have a non-trivially
colored crossing. Suppose at this crossing, the over-arc is colored by $b$ and the
two under-arcs are colored by $a$ and $c$. Then $2b=a+c$ (mod $p$) and $a, b, c$ are
distinct.

Let $f(x)=(b-a)^{-1}(x-a)$, where $x^{-1}$ is the inverse element of $x$ in $\mathbb{Z}_p$.
Then $f(a)=0$, $f(b)=1$ and $f(c)=\frac{c-a}{b-a}=\frac{2b-a-a}{b-a}=2$. Hence $f(C)$ is
a non-trivial $p$-coloring of $D$ containing colors 0, 1 and 2. \hfill{} $\Box$

\medskip

The fourth author of this article studied what kind of color set would never be a sufficient
set or minimal sufficient set. The following two theorems are useful.

\begin{theorem}\label{test1}\cite{Lopes}
Let $k$ be a positive integer and $L$ a link with non-zero determinant, admitting
non-trivial $(2k+1)$-colorings. If $S\subseteq \{0, 1, 2, \ldots, k\}$, then $S$
is not a $(2k+1)$-sufficient set of colors for $L$. Moreover, let $f$ be any
$(2k+1)$-coloring automorphism. Then any set $S'\subseteq \{f(0), f(1), . . . f(k)\}$
cannot be a $(2k+1)$-sufficient set of colors for L.
\end{theorem}

\begin{theorem}\label{test2}\cite{Lopes}
Let $k$ be a positive integer and $L$ a link with non-zero determinant, admitting
non-trivial $(2k+1)$-colorings. Suppose $S=\{c_1, \ldots, c_n\}$ is a $(2k+1)$-sufficient
set of colors for $L$. For each $c_i$ $(i=1, 2, \ldots, n)$, let $S_i$ be the set
of unordered pairs $\{c_i^1, c_i^2\}$ (with $c_i^1\neq c_i^2$ from $S$) such that
$2c_i=c_i^1+c_i^2 ~\mod~ 2k+1$.

$(1)$. There is at least one $i$ such that the expression $2c_i=c_i^1+c_i^2$ does not
make sense over the integers.

$(2)$. If there is $i\in \{1, 2, \ldots, n\}$ such that for all $j\in \{1, 2, \ldots, n\}\backslash \{i\}$,
$c_i$ is not in any of the unordered pairs in $S_j$, then $n> \text{mincol}_{2k+1} L$
and $S\backslash \{c_i\}$ is also a $(2k+1)$-sufficient set of colors for $L$.
\end{theorem}

\medskip

REMARK 2. Part $(1)$ of Theorem \ref{test2} is a generalization of Theorem \ref{test1}.

\medskip

Motivated by Theorem \ref{test2}, we define an associated edge-colored
simple graph for a color set.

\medskip

DEFINITION 6. ~Let $m$ be a positive integer. For a mod $m$ color set
$S=\{c_1, c_2, \ldots, c_n\}\subseteq \mathbb{Z}_{m}$, the associated edge-colored
(or briefly, colored) graph $G_c(S)$ is a simple graph constructed as follows.
Let $\{c_1, c_2, \ldots, c_n\}$ be the vertex set of $G_c(S)$. There is a red edge
between $c_i$ and $c_j$ if and only if there exists a positive integer $k\leq n$,
such that $2c_k=c_i+c_j$ over the integers. There is a blue edge between $c_i$ and
$c_j$ if and only if there exists an positive integer $k\leq n$, such that
$2c_k=c_i+c_j$ mod $m$ but $2c_k\neq c_i+c_j$ over the integers. We denote $G(S)$ the
underlying graph of $G_c(S)$, i.e., the graph obtained by replacing each colored
edge by a normal edge.

\medskip

\medskip

REMARK 3. Nakamura et al. defined a ``pallet graph'' for a $p$-coloring in \cite{NNS}. Their
graph and ours are similar in some respects. But our definition was made independently
and we have different motivations.

\medskip

Now we rewrite and extend Theorem \ref{test2} by using the associated graph.

\begin{theorem}\label{test3}
Let $k$ be a positive integer and $L$ a link with non-zero determinant, admitting
non-trivial $(2k+1)$-colorings. Suppose $S=\{c_1, \ldots, c_n\}$ is a $(2k+1)$-sufficient
set of colors for $L$. Let $G_c(S)$ ($G(S)$) be the associated colored (underlying)
graph.

$(1)$. $G_c(S)$ contains at least one blue edge.

$(2)$. If $c_i$ is an isolated vertex in $G(S)$, then $n> \text{mincol}_{2k+1} L$
and $S\backslash \{c_i\}$ is also a $(2k+1)$-sufficient set of colors for $L$.

$(3)$. Let $c(G(S))$ be the number of connected components of $G(S)$. Let $c(L)$
be the component number of $L$. If $c(G(S))> c(L)$, then $n> \text{mincol}_{2k+1} L$
and there exist $c(L)$ components of $G(S)$ such that the vertex set of these
$c(L)$ components is also a $(2k+1)$-sufficient set of colors for $L$.
\end{theorem}

\begin{corollary} (Corollary to both Theorem \ref{test2} and \ref{test3})\label{11}
$\{0, 1, 2, 3, 7\}$ is not a 11-minimal sufficient set of colors for any link
with non-zero determinant.
\end{corollary}

\medskip

\begin{figure}[htbp]
\centering
\includegraphics[width=10cm]{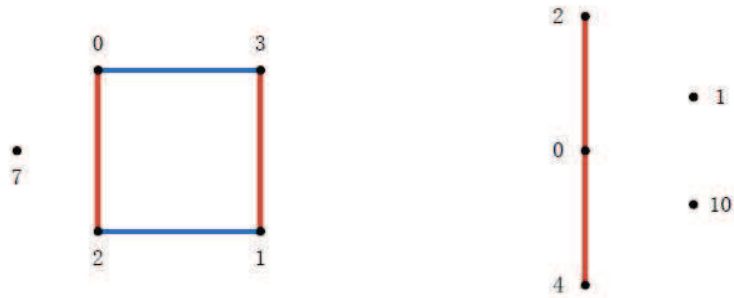}
\renewcommand{\figurename}{Fig.}
\caption{{\footnotesize The left is the colored graph associated with $\{0, 1, 2, 3, 7\}$,
the right is the colored graph associated with $\{0, 1, 2, 4, 10\}$.
}}\label{fig1}
\end{figure}

Proof. The graph associated with $\{0, 1, 2, 3, 7\}$ is shown in Fig. \ref{fig1}. According
to part (1) of Theorem \ref{test2} or Theorem \ref{test3}, $\{0, 1, 2, 3, 7\}$
is not a 11-minimal sufficient set of colors for any link with non-zero determinant. \hfill{} $\Box$

\medskip

\begin{corollary} (Corollary to both Theorem \ref{test2} and \ref{test3})\label{13}
$\{0, 1, 2, 4, 10\}$ is not a 13-minimal sufficient set of colors for any link
with non-zero determinant.
\end{corollary}

\medskip

Proof. The graph associated with $\{0, 1, 2, 4, 10\}$ is shown in Fig. \ref{fig1}. According
to either part (1) or part (2) of Theorem \ref{test2} or Theorem \ref{test3},
$\{0, 1, 2, 4, 10\}$ is not a 13-minimal sufficient set of colors for any link
with non-zero determinant. \hfill{} $\Box$

\medskip

The following proposition is obvious and we omit its proof here.

\begin{proposition}
Let $m$ be a positive integer. Let $S_1$ and $S_2$ be two mod $m$ color sets with the same
cardinality. If $S_1$ and $S_2$ are equivalent (mod $m$), then the associated underlying
graphs $G(S_1)$ and $G(S_2)$ are isomorphic.
\end{proposition}

\section{Minimal Sufficient Sets of Colors}
The fourth author of this article proved that unlike $p\leq 7$, there is no universal
$11$-minimal sufficient set of colors \cite{Lopes}. But it still remains unknown whether
there is a common $11$-minimal sufficient set of colors or not. We now give further results
on $11$-, $13$-, and $17$-sufficient or minimal sufficient sets of colors.

\subsection{Minimal Sufficient Sets of Colors with Cardinality $5$}

\begin{theorem}\label{p=11}
Let $L$ be an 11-colorable link with non-zero determinant. If a diagram of $L$ can be colored
by a color set of 5 colors, then the color set must be either $\{0, 1, 2, 3, 6\}$
or $\{0, 1, 2, 4, 7\}$ in the sense of equivalence of color sets induced by
coloring automorphism.
\end{theorem}

\medskip

Proof. It was shown in \cite{Lopes} that, in the sense of equivalence class of
color sets, knot $6_2$ has unique 11-minimal sufficient set $\{0, 1, 2, 3, 6\}$
and knot $7_2$ has unique 11-minimal sufficient set $\{0, 1, 2, 4,7\}$.
It was also shown in \cite{Lopes} that $\{0, 1, 2, 3, 6\}$ is not equivalent
to $\{0, 1, 2, 4, 7\}$. Recalling Lemma \ref{012}, we only need to consider
color sets of type $\{0, 1, 2, x, y\}$. Table 1 shows all instances of
color sets of type $\{0, 1, 2, x, y\}$. ``type 1'' (``type 2'') means it
is equivalent to $\{0, 1, 2, 3, 6\}$ ($\{0, 1, 2, 4, 7\}$). ``N, Th \ref{test1}''
(``N, Co \ref{11}'') means there is no link which can be 11-colored by it due
to Theorem \ref{test1} (Corollary \ref{11}). We pick three color sets as
examples of how to read Table 1.

(1) $\{0, 1, 2, 3, 8\}$. The coloring automorphism $f_{10, 3}(x)=10x+3$
transforms the color set $\{0, 1, 2, 3, 8\}$ into $\{0, 1, 2, 3, 6\}$,
so $\{0, 1, 2, 3, 8\}$ is equivalent to $\{0, 1, 2, 3, 6\}$.

(2) $\{0, 1, 2, 3, 9\}$. The coloring automorphism $f_{1, 2}(x)=x+2$
transforms the color set $\{0, 1, 2, 3, 9\}$ into $\{0, 2, 3, 4, 5\}$, a
subset of $\{0, 1, \ldots, 5\}$. According to Theorem \ref{test1},
$\{0, 1, 2, 3, 9\}$ is not $11$-sufficient set of colors for any link
with non-zero determinant.

(3) $\{0, 1, 2, 4, 9\}$. The coloring automorphism $f_{6, 1}(x)=6x+1$
transforms the color set $\{0, 1, 2, 4, 9\}$ into $\{0, 1, 2, 3, 7\}$.
According to Corollary \ref{11}, $\{0, 1, 2, 3, 7\}$ is not a $11$-minimal
sufficient set of colors for any link with non-zero determinant. Hence
$\{0, 1, 2, 4, 9\}$ is not a $11$-minimal sufficient set of colors for
any link with non-zero determinant either. \hfill{} $\Box$

\medskip

\begin{table}[ht]
\center{Table 1. 5 color sets mod 11.}
\vskip 5 true pt
{\begin{tabular}{|l|r|c||l|r|c|}
 \hline
0,1,2,3,4  &                          & N, Th \ref{test1} & 0,1,2,5,7  & ($\times$2+1)=0,1,3,4,5   & N, Th \ref{test1} \\ \hline
0,1,2,3,5  &                          & N, Th \ref{test1} & 0,1,2,5,8  & ($\times$7)=0,1,2,3,7     & N, Co \ref{11}    \\ \hline
0,1,2,3,6  &                          & type 1            & 0,1,2,5,9  & ($\times$8+6)=0,1,2,3,6   & type 1            \\ \hline
0,1,2,3,7  &                          & N, Co \ref{11}    & 0,1,2,5,10 & (+1)=0,1,2,3,6            & type 1            \\ \hline
0,1,2,3,8  & ($\times$10+3)=0,1,2,3,6 & type 1            & 0,1,2,6,7  & ($\times$2)=0,1,2,3,4     & N, Th \ref{test1} \\ \hline
0,1,2,3,9  & (+2)=0,2,3,4,5           & N, Th \ref{test1} & 0,1,2,6,8  & ($\times$2)=0,1,2,4,5     & N, Th \ref{test1} \\ \hline
0,1,2,3,10 & (+1)=0,1,2,3,4           & N, Th \ref{test1} & 0,1,2,6,9  & ($\times$2)=0,1,2,4,7     & type 2            \\ \hline
0,1,2,4,5  &                          & N, Th \ref{test1} & 0,1,2,6,10 & (+1)=0,1,2,3,7            & N, Co \ref{11}    \\ \hline
0,1,2,4,6  & ($\times$6)=0,1,2,3,6    & type 1            & 0,1,2,7,8  & ($\times$2)=0,2,3,4,5     & N, Th \ref{test1} \\ \hline
0,1,2,4,7  &                          & type 2            & 0,1,2,7,9  & ($\times$5+1)=0,1,2,3,6   & type 1            \\ \hline
0,1,2,4,8  & ($\times$3)=0,1,2,3,6    & type 1            & 0,1,2,7,10 & ($\times$10+2)=0,1,2,3,6  & type 1            \\ \hline
0,1,2,4,9  & ($\times$6+1)=0,1,2,3,7  & N, Co \ref{11}    & 0,1,2,8,9  & (+3)=0,1,3,4,5            & N, Th \ref{test1} \\ \hline
0,1,2,4,10 & (+1)=0,1,2,3,5           & N, Th \ref{test1} & 0,1,2,8,10 & (+3)=0,2,3,4,5            & N, Th \ref{test1} \\ \hline
0,1,2,5,6  & ($\times$2+1)=0,1,2,3,5  & N, Th \ref{test1} & 0,1,2,9,10 & (+2)=0,1,2,3,4            & N, Th \ref{test1} \\ \hline
\end{tabular}}
\end{table}

\begin{theorem}\label{p=13}
Let $L$ be a $13$-colorable link with non-zero determinant. If a diagram of $L$ can be
colored by 5 colors, then it can be only colored by $\{0, 1, 2, 4, 7\}$ in
the sense of equivalence. Specifically, if $\text{mincol}_{13} L=5$ holds
for all 13-colorable links with non-zero determinant, then $\{0, 1, 2, 4, 7\}$
is the only common 13-minimal sufficient set of colors in the sense of
equivalence of color sets induced by coloring automorphism.
\end{theorem}

\medskip

Proof. It was shown in \cite{Lopes} that, in the sense of equivalence class
of color sets, knots $6_3$, $7_3$ and $10_{154}$ has unique 13-minimal
sufficient set $\{0, 1, 2, 4, 7\}$. Recalling Lemma \ref{012}, we only
need to consider color sets of type $\{0, 1, 2, x, y\}$. Table 2 shows
circumstances of all color sets of type $\{0, 1, 2, x, y\}$. ``Y'' means
it is equivalent to $\{0, 1, 2, 4, 7\}$. ``N, Th \ref{test1}'' (``N, Co \ref{13}'')
means there is no link can be 11-colored by it due to Theorem \ref{test1}
(Corollary \ref{13}). \hfill{} $\Box$

\medskip

\begin{table}[ht]
\center{Table 2. 5 color sets mod 13.}
\vskip 5 true pt
{\begin{tabular}{|l|r|c||l|r|c|}
\hline
0,1,2,3,4  &                            & N, Th \ref{test1} & 0,1,2,5,12  & (+1)=0,1,2,3,6            & N, Th \ref{test1} \\ \hline
0,1,2,3,5  &                            & N, Th \ref{test1} & 0,1,2,6,7   & ($\times$2+1)=0,1,2,3,5   & N, Th \ref{test1} \\ \hline
0,1,2,3,6  &                            & N, Th \ref{test1} & 0,1,2,6,8   & ($\times$2+1)=0,1,3,4,5   & N, Th \ref{test1} \\ \hline
0,1,2,3,7  & ($\times$2)=0,1,2,4,6      & N, Th \ref{test1} & 0,1,2,6,9   & ($\times$3)=0,1,3,5,6     & N, Th \ref{test1} \\ \hline
0,1,2,3,8  & ($\times$2)=0,2,3,4,6      & N, Th \ref{test1} & 0,1,2,6,10  & ($\times$3)=0,3,4,5,6     & N, Th \ref{test1} \\ \hline
0,1,2,3,9  & ($\times$2)=0,2,4,5,6      & N, Th \ref{test1} & 0,1,2,6,11  & ($\times$10+7)=0,1,2,4,7  & Y                 \\ \hline
0,1,2,3,10 & (+3)=0,3,4,5,6             & N, Th \ref{test1} & 0,1,2,6,12  & ($\times$2+2)=0,1,2,4,6   & N, Th \ref{test1} \\ \hline
0,1,2,3,11 & (+2)=0,2,3,4,5             & N, Th \ref{test1} & 0,1,2,7,8   & ($\times$2)=0,1,2,3,4     & N, Th \ref{test1} \\ \hline
0,1,2,3,12 & (+1)=0,1,2,3,4             & N, Th \ref{test1} & 0,1,2,7,9   & ($\times$2)=0,1,2,4,5     & N, Th \ref{test1} \\ \hline
0,1,2,4,5  &                            & N, Th \ref{test1} & 0,1,2,7,10  & ($\times$2)=0,1,2,4,7     & Y                 \\ \hline
0,1,2,4,6  &                            & N, Th \ref{test1} & 0,1,2,7,11  & ($\times$6+1)=0,1,2,4,7   & Y                 \\ \hline
0,1,2,4,7  &                            & Y                 & 0,1,2,7,12  & ($\times$2+2)=0,2,3,4,6   & N, Th \ref{test1} \\ \hline
0,1,2,4,8  & ($\times$7)=0,1,2,4,7      & Y                 & 0,1,2,8,9   & ($\times$2)=0,2,3,4,5     & N, Th \ref{test1} \\ \hline
0,1,2,4,9  & ($\times$3+1)=0,1,2,4,7    & Y                 & 0,1,2,8,10  & ($\times$11+4)=0,1,2,4,10 & N, Co \ref{13}    \\ \hline
0,1,2,4,10 &                            & N, Co \ref{13}    & 0,1,2,8,11  & ($\times$12+2)=0,1,2,4,7  & Y                 \\ \hline
0,1,2,4,11 & (+2)=0,2,3,4,6             & N, Th \ref{test1} & 0,1,2,8,12  & ($\times$2+2)=0,2,4,5,6   & N, Th \ref{test1} \\ \hline
0,1,2,4,12 & (+1)=0,1,2,3,5             & N, Th \ref{test1} & 0,1,2,9,10  & (+4)=0,1,4,5,5            & N, Th \ref{test1} \\ \hline
0,1,2,5,6  &                            & N, Th \ref{test1} & 0,1,2,9,11  & (+4)=0,2,4,5,6            & N, Th \ref{test1} \\ \hline
0,1,2,5,7  & ($\times$2)=0,1,2,4,10     & N, Co \ref{13}    & 0,1,2,9,12  & (+4)=0,3,4,5,6            & N, Th \ref{test1} \\ \hline
0,1,2,5,8  & ($\times$11+4)=0,1,2,4,7   & Y                 & 0,1,2,10,11 & (+3)=0,1,3,4,5            & N, Th \ref{test1} \\ \hline
0,1,2,5,9  & ($\times$3)=0,1,2,3,6      & N, Th \ref{test1} & 0,1,2,10,12 & (+3)=0,2,3,4,5            & N, Th \ref{test1} \\ \hline
0,1,2,5,10 & ($\times$3)=0,2,3,4,6      & N, Th \ref{test1} & 0,1,2,11,12 & (+2)=0,1,2,3,4            & N, Th \ref{test1} \\ \hline
0,1,2,5,11 & ($\times$12+2)=0,1,2,4,10  & N, Co \ref{13}    &             &                           &                   \\ \hline
\end{tabular}}
\end{table}

\subsection{Possible Minimal Sufficient Sets of Colors with Cardinality $6$}

Now we determine which color sets with cardinality $6$ may be minimal sufficient sets of colors
for primes $p=11$, $13$, and $17$. Our strategy is as follows.

{\bf Step 1.} List all subsets of $\{0, 1, \ldots, p\}$ with cardinality $6$ and
containing $0, 1$, and $2$.

{\bf Step 2.} Classify these color sets into equivalence classes (recall
Definition 5).

{\bf Step 3.} Use Theorem \ref{test1} and Theorem \ref{test2} to check the color
sets in order. If $S$ can not be a $p$-minimal sufficient set of colors
for any link with non-zero determinant, then delete all those color sets
equivalent to $S$. We call the remaining ones possible $p$-sufficient
sets of colors with cardinality $6$.

Our results follow. We use a C program to achieve Step 1 and 2.

For $p=11$, there are $\binom{8}{3}=56$ color sets with cardinality $6$ and
containing $0, 1$, and $2$. They can be classified into $6$ equivalence classes.
Among them, there are $4$ possible $11$-minimal sufficient sets of colors up
to the equivalence relation: $\{0, 1, 2, 3, 4, 6\}$, $\{0, 1, 2, 3, 4, 7\}$,
$\{0, 1, 2, 3, 5, 6\}$, and $\{0, 1, 2, 3, 5, 9\}$.

\begin{proposition}
To the extent encompassed by Theorems \ref{test1} and \ref{test2},
there are no obstructions for an $11$-colorable link $L$ with non-zero
determinant to satisfy $\text{mincol}_{11} L = 5$.
\end{proposition}

\medskip

Proof. By \cite{CJZ}, $\text{mincol}_{11} L$ is either $5$ (in which case
the proof is concluded) or $6$. In the latter instance, there is a
diagram of $L$ equipped with a non-trivial $11$-coloring using $6$
colors either from $\{0, 1, 2, 3, 4, 6\}$, or from $\{0, 1, 2, 3, 4, 7\}$,
or from $\{0, 1, 2, 3, 5, 6\}$, or from $\{0, 1, 2, 3, 5, 9\}$, up
to the equivalence relation. We now prove that it is always possible
to remove one color from any of these sets in order to obtain
$\{0, 1, 2, 3, 6\}$ or $\{0, 1, 2, 4, 7\}$. Using the tests described in
Theorems \ref{test1} and/or \ref{test2} along with the help of Table 1,
we note that $3$ or $4$ can be removed from the set $\{0, 1, 2, 3, 4, 6\}$
giving rise to a set equivalent to $\{0, 1, 2, 3, 6\}$;
$0$ or $1$ ($3$) can be removed from the set $\{0, 1, 2, 3, 4, 7\}$ giving
rise to a set equivalent to $\{0, 1, 2, 3, 6\}$ ($\{0, 1, 2, 4, 7\}$);
$5$ can be removed from the set $\{0, 1, 2, 3, 5, 6\}$ giving rise to
$\{0, 1, 2, 3, 6\}$; $0$ or $3$ can be removed from the set $\{0, 1, 2, 3, 5, 9\}$
giving rise to a set equivalent to $\{0, 1, 2, 3, 6\}$. End of proof. \hfill{} $\Box$

\medskip

For $p=13$, there are $\binom{10}{3}=120$ color sets with cardinality 6 and
containing 0, 1, and 2. They can be classified into 14 equivalence classes.
Among them, there are 8 possible 13-minimal sufficient sets of colors up to
the equivalence relation: \{0, 1, 2, 3, 4, 7\}, \{0, 1, 2, 3, 4, 8\},
\{0, 1, 2, 3, 5, 7\}, \{0, 1, 2, 3, 5, 8\}, \{0, 1, 2, 3, 5, 9\},
\{0, 1, 2, 3, 5, 11\}, \{0, 1, 2, 3, 6, 10\}, and \{0, 1, 2, 4, 5, 8\}.

For $p=17$, there are $\binom{14}{3}=364$ color sets with cardinality 6 and
containing 0, 1, and 2. They can be classified into 49 equivalence classes.
Among them, there are 9 possible 17-minimal sufficient sets of colors up to
the equivalence relation: \{0, 1, 2, 3, 5, 9\}, \{0, 1, 2, 3, 5, 10\},
\{0, 1, 2, 3, 5, 12\}, \{0, 1, 2, 3, 6, 9\}, \{0, 1, 2, 3, 6, 10\},
\{0, 1, 2, 3, 6, 11\}, \{0, 1, 2, 3, 6, 13\}, \{0, 1, 2, 3, 6, 14\},
and \{0, 1, 2, 3, 7, 10\}.

\section{Minimum Number of Colors}
In this section, we study the lower bound of the minimum number of colors.
First we give a short proof of $\text{mincol}_p L\geq 5$ for any
$p$-colorable link $L$ with $\det L\neq 0$ and prime $p\geq 11$. Then
we show that we can go further by using a similar approach.

\begin{theorem}\label{lowerbound5}\cite{LM}
Let $L$ be a link with non-zero determinant. If there is a prime $p\geq 11$
such that $L$ admits non-trivial $p$-colorings, then $\text{mincol}_p L\geq 5$.
\end{theorem}

\medskip

Proof. It is easy to prove the mod $p$ minimum number of colors is at least $4$ for
any link with non-zero determinant and any prime $p\geq 5$. See Lemma 2.1 in
\cite{Satoh} for example. So we only need to prove $L$ cannot be colored by
4 colors. Recalling Lemma \ref{012}, it is enough to consider color set $S$ of
type $\{0, 1, 2, a\}$, where $3\leq a \leq p-1$. If $L$ can be colored by $S$,
then $\frac{p+1}{2} \leq a \leq p-1$ due to Theorem \ref{test1}. The coloring
automorphism $f(x)=x+(p-a)$ transforms $S$ into $\{0, p-a, p-a+1, p-a+2\}$.
Hence $p-a+2\geq \frac{p-1}{2}+1$ due to Theorem \ref{test1}, thus $a\leq \frac{p+3}{2}$.
So $a\in \{\frac{p+1}{2}, \frac{p+3}{2}\}$. We use the same notations as in
Theorem \ref{test2}. Let $c_1=0$, $c_2=1$, $c_3=2$, $c_4=a$. Since
$a\in \{\frac{p+1}{2}, \frac{p+3}{2}\}$, we have
$6\leq \frac{p+1}{2} \leq c_4 \leq \frac{p+3}{2} = p-\frac{p-3}{2} \leq p-4$.
It is easy to check $S_1=\emptyset$, $S_2=\{\{0, 2\}\}$ and $S_3=\emptyset$.
So $c_4$ is not in any of the unordered pairs in $S_j$, $j\in \{1, 2, 3\}$.
According to part $(2)$ of Theorem \ref{test2}, $4> \text{mincol}_p L$, which
never happens. \hfill{} $\Box$

\medskip

\begin{theorem}\label{lowerbound6}
Let $L$ be a link with non-zero determinant. If there is a prime $p\geq 17$
such that $L$ admits non-trivial $p$-colorings, then $\text{mincol}_p L\geq 6$.
\end{theorem}

\medskip

Proof. Recalling part $(5)$ of Theorem \ref{mincol}, we only need to prove $L$ cannot
be colored by 5 colors. Suppose $L$ can be colored by a color set $S$ with
$|S|=5$. Then $S$ is a $p$-minimal sufficient set of colors for $L$.

\medskip

{\bf Claim.} Suppose $S'=\{x_1, x_2, x_3, x_4\}\subseteq S$ with $x_1< x_2< x_3< x_4$.
Then $x_4-x_1\geq 5$.

\medskip

Proof of the Claim. We prove the Claim by contradiction. If there is an $S'$ such that $x_4-x_1\leq 4$,
then the coloring automorphism $f_{1, p-x_1}(x)=x+p-x_1$ transforms $S$ into a
set $S''=\{0, y_1, y_2, y_3, y_4\}$, where $y_1< y_2< y_3\leq 4$ and
$y_3< y_4$. According to Theorem \ref{test1}, $y_4\geq \frac{p+1}{2}$ and
$p-y_4+4\geq \frac{p+1}{2}$, which yields $9\leq y_4 \leq p-5$.
So $y_4$ is not in any of the unordered pairs in any $S_j$ (which means $y_4$
cannot be the color of any under-arc at a polychromatic crossing). Recalling
Theorem \ref{test2}, $S''$ (so that $S$) is not a $p$-minimal sufficient set
of colors for $L$ and $\text{mincol}_p L\leq 4$, which is impossible. The
proof is complete.

\medskip

Recalling Lemma \ref{012}, it is enough to consider color set $S$ of type
$\{0, 1, 2, a, b\}$, where $3\leq a \leq b \leq p-1$. Then
$\frac{p+1}{2} \leq b \leq p-1$ due to Theorem \ref{test1}.

The Claim indicates that $a\geq 5$ immediately. It also indicates that
$b \leq p-3$, otherwise the coloring automorphism $f_{1, p-b}(x)=x+p-b$
will transform $\{0, 1, 2, b\}$ into $\{0, p-b, p-b+1, p-b+2\}$ where
$p-b+2\leq 4$. Hence $a\geq 5$ and $\frac{p+1}{2}\leq b \leq p-3$. So
there is no crossing with the over-arc colored by one color from $\{0, 1, 2\}$,
whose under-arcs are colored with one color from $\{0, 1, 2\}$, and the other from $\{a, b\}$.

We keep all the notations used in Theorem \ref{test2}. Let $c_1=0$, $c_2=1$,
$c_3=2$, $c_4=a$, $c_5=b$. Now we divide all the possibilities into the
following two cases.

\noindent{\bf Case 1.} $a+b\notin \{0, 2, 4\}\mod p$.

In this case, both $a$ and $b$ are not in any of the unordered pairs in
$S_1$, $S_2$ and $S_3$. So $a$ must be in one of the unordered pairs in
$S_5$ and $b$ must be in one of the unordered pairs in $S_4$ according
to Theorem \ref{test2}, i.e., $\exists \alpha, \beta \in \{0, 1, 2\}$($\alpha$ and $\beta$ are not
necessarily different), such that
\begin{equation}\label{1}
\left\{ \begin{array}{l}
2b=a+\alpha+p \\
2a=b+\beta
\end{array}. \right.
\end{equation}
Hence
\begin{equation}\label{2}
\left\{ \begin{array}{l}
a=\frac{p+\alpha+2\beta}{3} \\
b=\frac{2p+2\alpha+\beta}{3}
\end{array}. \right.
\end{equation}

The system of equations (\ref{1}) is over integers, so are
(\ref{3}) and (\ref{5}) below. Since
$\alpha, \beta \in \{0, 1, 2\}$, the coloring automorphism $f_{3, 0}(x)=3x$
transforms $S$ into $\{0, 3, 6, \alpha+2\beta, 2\alpha+\beta\}\subseteq \{0, 1, \ldots, \frac{p-1}{2}\}$.
So $S$ is not a $p$-sufficient set of colors for $L$.

\noindent{\bf Case 2.} $a+b=\gamma \in \{0, 2, 4\}\mod p$.

It is easy to see $c_2=1$ is not in any of the unordered pairs in $S_1$
and $S_3$. So $1$ must be in one of the unordered pairs in $S_4$ or $S_5$,
i.e., either $2a=1+b$ or $2b=1+a+p$ (over integers). So there are two subcases.

\noindent{\bf Case 2.1.} $2a=1+b$.

In this subcase,
\begin{equation}\label{3}
\left\{ \begin{array}{l}
2a=1+b \\
a+b=\gamma+p
\end{array}. \right.
\end{equation}
Hence
\begin{equation}\label{4}
\left\{ \begin{array}{l}
a=\frac{p+\gamma+1}{3} \\
b=\frac{2p+2\gamma-1}{3}
\end{array}. \right.
\end{equation}
Since $\gamma \in \{0, 2, 4\}$, the coloring automorphism $f_{3, 1}(x)=3x+1$
transforms $S$ into $\{1, 4, 7, \gamma+2, 2\gamma\}\subseteq \{0, 1, \ldots, \frac{p-1}{2}\}$.
So $S$ is not a $p$-sufficient set of colors for $L$.

\noindent{\bf Case 2.2.} $2b=1+a+p$.

In this subcase,
\begin{equation}\label{5}
\left\{ \begin{array}{l}
2b=1+a+p \\
a+b=\gamma+p
\end{array}. \right.
\end{equation}
Hence
\begin{equation}\label{6}
\left\{ \begin{array}{l}
a=\frac{p+2\gamma-1}{3} \\
b=\frac{2p+\gamma+1}{3}
\end{array}. \right.
\end{equation}
Since $\gamma \in \{0, 2, 4\}$, the coloring automorphism $f_{3, 1}(x)=3x+1$
transforms $S$ into $\{1, 4, 7, 2\gamma, \gamma+2\}\subseteq \{0, 1, \ldots, \frac{p-1}{2}\}$.
So $S$ is not a $p$-sufficient set of colors for $L$.

Therefore, $S$ is not a $p$-minimal sufficient set of colors for $L$, a contradiction. \hfill{} $\Box$

\medskip

As we pointed out before, for any split link and any modulus, the
minimum number of colors is $2$. So Theorem \ref{lowerbound5},
Theorem \ref{lowerbound6} and Theorem \ref{NNS} cannot be extended
to all links. But can we use non-split links instead of links with
non-zero determinant in these theorems? The following lemma gives
a negative answer.

\begin{lemma}
For any modulus $m\geq 5$, there are infinitely many non-split links
with zero determinant having minimum number of colors at most $4$.
\end{lemma}

\medskip

\begin{figure}[htbp]
\centering
\includegraphics[width=5cm]{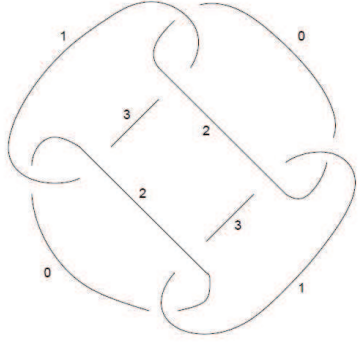}
\renewcommand{\figurename}{Fig.}
\caption{{\footnotesize A non-trivial coloring of $L8n8$ with $4$ colors.
}}\label{L8n8}
\end{figure}

Proof. Figure \ref{L8n8} shows that for any modulus $m\geq 5$, the non-split link
$L8n8$ in the Thistlethwaite link table has a non-trivial coloring with
color set $\{0, 1, 2, 3\}$. Note that $L8n8$ is $(2, -2, 2, -2)$-pretzel
link. It is easy to see that for any $n$,
$(\underbrace{2, -2, \ldots, 2, -2}_{2n~\text{strands}})$-pretzel link
has the same property.

For pretzel link $P(p_1, p_2, \ldots, p_n)$, the determinant is
$\Big|\sum_{j=1}^n \frac{\prod_{i=1}^n p_i}{p_j}\Big|$ \cite{DFKLS, GZ}.
Hence $\det (P(2, -2, \ldots, 2, -2))=0$. This completes the proof. \hfill{} $\Box$

\medskip

\section{A Question Raised By Nakamura, Nakanishi, and Satoh}
We recall that Theorem \ref{NNS} which states that
$\text{mincol}_p (K) \geq \lfloor \log_2 p \rfloor +2$ is proved in
Nakamura et al.' \cite{NNS}. Since the crossing number of knot $K$,
$c(K)$, satisfies $c(K)\geq \text{mincol}_p (K)$, for any $p$-colorable
knot $K$, these authors wonder if the equality $c(K)= \lfloor \log_2 p \rfloor +2$
only holds for the trefoil and the figure-eight knots, see (iii) in Remark
3.3 on page 96 of \cite{NNS}. Here we settle this matter with Theorem \ref{equality}.

\begin{theorem}\label{equality}
Let $p$ be an odd prime. Let $K$ be a $p$-colorable classical knot. Then the equality in
$c(K)\geq \lfloor \log_2 p \rfloor +2$ only holds for the trefoil knot ($p=3$)
and the figure-eight knot ($p=5$).
\end{theorem}

Let $D$ be a link diagram. Let
$$d_n^\infty := \max \{\det(D)\ |\ D\ \text{is a link diagram of $n$ crossings}\}.$$
In \cite{Sto}, Stoimenow showed
\begin{eqnarray}
d_n^\infty\leq d_{n-1}^\infty + d_{n-2}^\infty + d_{n-3}^\infty\ \ (n>2), \label{f1}
\end{eqnarray}
and then proved the following theorem.

\begin{theorem}\cite{Sto}\label{Sto}
Let $\delta\approx 1.83929$ be the real positive root of $x^3-x^2-x-1=0$.
There exists a constant $C>0$ such that for any link diagram $D$ of $c(D)$
crossings
$$\det(D)\leq C\cdot \delta^{c(D)}.$$
\end{theorem}

Stoimenow also pointed out that $C=1$ is always valid.

\medskip

\noindent{Proof of Theorem \ref{equality}.}

Let $\tilde{D}$ be a minimal diagram of $K$. Since $K$ is a $p$-colorable
knot, we have $p \mid \det(K)$ and $\det K>0$. By Theorem \ref{Sto},
$$\log_2 p \leq \log_2 \det(K)=\log_2 \det(\tilde{D})\leq c(\tilde{D})\log_2 \delta <0.87915\cdot c(K).$$
It is easy to see, for $c(K)\geq 17$,
$$c(K)>0.87915\cdot c(K)+2>\log_2 p+2>\lfloor \log_2 p \rfloor +2.$$

Table 3 shows the numerical results of $d_n^\infty$ and
$\lfloor \log_2 d_n^\infty \rfloor +2$ for $3\leq n\leq 16$. The first
four values of $d_n^\infty$ ($n\geq 3$) appeared in \cite{Sto} and other
values are estimated by formula (\ref{f1}).

\begin{table}[ht]
\center{Table 3. $d_n^\infty$ and $\lfloor \log_2 d_n^\infty \rfloor +2$ for $3\leq n\leq 16$.}
\vskip 5 true pt
{\begin{tabular}{|c|c|c|c|c|c|c|c|}
\hline
$n$                                     &    $3$     &    $4$     &    $5$     &    $6$      &    $7$      &    $8$      &     $9$     \\  \hline
$d_n^\infty$                            &    $3$     &    $5$     &    $8$     &    $16$     & $\leq 29$   & $\leq 53$   & $\leq 98$   \\  \hline
$\lfloor \log_2 d_n^\infty \rfloor +2$  &    $3$     &    $4$     &    $5$     &    $6$      & $\leq 6$    & $\leq 7$    & $\leq 8$    \\  \hline
$n$                                     &    $10$    &    $11$    &    $12$    &    $13$     &    $14$     &    $15$     &     $16$    \\  \hline
$d_n^\infty$                            & $\leq 180$ & $\leq 331$ & $\leq 609$ & $\leq 1120$ & $\leq 2060$ & $\leq 3789$ & $\leq 6969$ \\  \hline
$\lfloor \log_2 d_n^\infty \rfloor +2$  & $\leq 9$   & $\leq 10$  & $\leq 11$  & $\leq 12$   & $\leq 13$   & $\leq 13$   & $\leq 14$   \\  \hline
\end{tabular}}
\end{table}

Hence, for any knot $K$ with crossing number between $7$ and $16$, we obtain
$$c(K)>\lfloor \log_2 d_n^\infty \rfloor +2\geq \lfloor \log_2 \det(K) \rfloor +2\geq \lfloor \log_2 p \rfloor +2.$$

For any knot $K$ with crossing number $5$ or $6$, it is easy to check that
$c(K)>\lfloor \log_2 p \rfloor +2$. The proof is complete. {\hfill$\Box$}

\vskip 20 true pt

\noindent ACKNOWLEDGEMENTS. ~X. Jin and J. Ge acknowledge support from the National Natural Science
Foundation of China (No. 11271307). P. Lopes acknowledges partial
funding from FCT (Portugal) through projects PEst-OE/EEI/LA0009/2013,
and EXCL/MAT-GEO/0222/2012 (``Geometry and Mathematical Physics'').
L. Zhang acknowledges support from the National Natural Science Foundation
of China (Nos. 11171279 and 11471273). The authors thank Prof. Satoh for
kindly sending us \cite {NNS}. J. Ge thanks Jingkang Zhou for helping
writing the program.

\vskip 20 true pt

\noindent Jun Ge \\
School of Mathematical Sciences, \\
Xiamen University, \\
Xiamen, Fujian 361005, P. R. China \\
Email: mathsgejun@163.com

\bigskip

\noindent Xian'an Jin \\
School of Mathematical Sciences, \\
Xiamen University, \\
Xiamen, Fujian 361005, P. R. China \\
Email: xajin@xmu.edu.cn

\bigskip

\noindent Louis H. Kauffman \\
Department of Mathematics, Statistics and Computer Science, \\
University of Illinois at Chicago, \\
851 S. Morgan St., Chicago IL 60607-7045, USA \\
Email: kauffman@uic.edu

\bigskip

\noindent Pedro Lopes \\
Center for Mathematical Analysis, Geometry, and Dynamical Systems, \\
Department of Mathematics, \\
Instituto Superior T\'{e}cnico, \\
Universidade de Lisboa,\\
1049--001 Lisbon, Portugal \\
Email: pelopes@math.tecnico.ulisboa.pt

\bigskip

\noindent Lianzhu Zhang \\
School of Mathematical Sciences, \\
Xiamen University, \\
Xiamen, Fujian 361005, P. R. China \\
Email: zhanglz@xmu.edu.cn

\end{document}